\documentclass[11pt]{article}

\usepackage{amssymb}
\usepackage{amsmath}
\usepackage{theorem}
\usepackage{epsfig}
\usepackage{verbatim}
\usepackage{graphicx}

\newtheorem{thm}{Theorem}[section]

\newtheorem{lemma}[thm]{Lemma}

\newcommand{\R}{\Bbb{R}}

\newcommand{\D}{\displaystyle}

\newcommand{\grad}{\nabla}

\newcommand{\dt}{\frac{d}{dt}}

\newcommand{\dxu}{\partial_{x_1}}
\newcommand{\dxd}{\partial_{x_2}}
\newcommand{\dxt}{\partial_{x_3}}

\newcommand{\al}{\alpha}

\begin{document}

\author{Diego C\'ordoba and Francisco Gancedo}
\title{Absence of squirt singularities for the multi-phase\\ Muskat problem}
\date{November 20, 2009}

\maketitle

\begin{abstract}
In this paper we study the evolution of multiple fluids with different constant densities in porous media. This physical scenario is known as the Muskat and the (multi-phase) Hele-Shaw problems. In this context we prove that the fluids do not develop squirt singularities.
\end{abstract}

\maketitle

\section{Introduction}
We consider the dynamics of the two interphases between three incompressible and immiscible fluids in a porous medium without surface tension. The free boundaries are prescribed by the jump of densities between the fluids. Although the present paper is devoted to the evolution of two interphases in $\R^3$, the same approach can be performed when additional phases are considered in $\R^2$ or $\R^3$. The governing equations for  the 2D incompressible porous media are identical with those modeling the dynamics in a Hele-Shaw cell (see \cite{S-T}).

The precise formulation of this problem is as follows \cite{bear}: the scalar density $\rho =\rho (x,t)$ of the fluid is convected by the incompressible velocity flow $u = (u_1(x,t), u_2(x,t), u_3(x,t))$ which satisfies Darcy's law i.e.
\begin{equation}\label{efm}
\begin{array}{ll}
\rho_t + u\cdot\nabla\rho =0& \quad\quad \text{(Conservation of mass)}\\
\quad\nabla\cdot u =0&             \quad\quad \text{(Incompressibility)}\\
u=-\grad p-(0,0, \rho)&       \quad\quad \text{(Darcy's law)}\\
\end{array}
\end{equation}
where the scalar $p=p(x,t)$ is the pressure and the acceleration due to gravity is taken equal to one to simplify the notation.

Darcy's law yields the velocity in terms of the density by means of singular integral operators as follows:
\begin{eqnarray}\label{sio}
    u(x,t)= PV \int_{\R^3} K(x-y)\,\rho(y,t)dy\, - \frac{2}{3}\left(0, 0,\rho(x)\right),\quad x\in \R^3,
\end{eqnarray}
where the kernel $K$ is given by
$$
K(x)=\frac{1}{4\pi}\left(3\frac{ x_1
x_3}{|x|^5}, 3\frac{ x_2
x_3}{|x|^5}, \frac{2x^2_3 - x_1^2-x_2^2}{|x|^5}\right).
$$
The integral operator is defined in the Fourier side by
\begin{eqnarray*}
\widehat{u}(\xi)=
(\frac{\xi_1\xi_3}{|\xi|^2},\frac{\xi_2\xi_3}{|\xi|^2},-\frac{\xi_1^2 + \xi_2^2}{|\xi|^2})\widehat{\rho}(\xi)
\end{eqnarray*}
that shows that the velocity and the density are at the same level in terms of regularity.

The fluid is characterized by three different constant values of density $\rho^1 < \rho^2 < \rho^3$
\begin{equation*}
 \rho(x_1,x_2,x_3,t)=\left\{\begin{array}{cl}
                    \rho^1\quad\mbox{in}&\Omega^1=\{x_3>f(x_1,x_2,t)\},\\
                    \rho^2\quad\mbox{in}&\Omega^2=\{f(x_1,x_2,t)>x_3>g(x_1,x_2,t)\},\\
                    \rho^3\quad\mbox{in}&\Omega^3=\{g(x_1,x_2,t)>x_3\},

                 \end{array}\right.
\end{equation*}
where $f(x_1,x_2,t)> g(x_1,x_2,t)$.  The moving surfaces
$$
S_f (x_1,x_2,t ) = \{(x_1 , x_2 , x_3 ) \in \R^3: x_3 = f ( x_1 , x_2 , t ) \}
$$
$$
S_g (x_1,x_2, t ) = \{(x_1 , x_2 , x_3 ) \in \R^3: x_3 = g ( x_1 , x_2 , t ) \}
$$
have the property (see below in section 2) that they can be  parameterized  as a graph for all time. Since the flow is incompressible the velocity of each interphase is continuous in the normal direction to the moving surface. Moreover, from the formulation it implies that the pressures are equal across the surfaces. The system is highly non-local in the sense that the equation for the moving surfaces involves singular integral operators and they are coupled together. Within the formulation we can recover the dynamics of a single interphase by taking $\rho^1 = \rho^2$ or $\rho^2 = \rho^3$, which has been shown to be well-possed in the stable scenario with a maximum principle (see \cite{DY} and \cite{DY2}). This case is known as the Muskat problem \cite{Muskat} (in 2D is also known as the two-phase Hele-Shaw)  which has been broadly studied in \cite{Peter}, \cite{ES}, \cite{DPS}, \cite{SCH}, \cite{DY}, \cite{DY2} and reference therein.

The aim of this paper is first to show that for multiple interphases, in the stable case ($\rho^1 < \rho^2 < \rho^3$ ), the system is well possed in a chain of Sobolev spaces. In the unstable case ($\rho^1 > \rho^2$ or $\rho^2 > \rho^3$ ) there is Rayleigh-Taylor instability \cite{DY} and the system is ill-possed. Secondly we rule out a squirt singularity in the three phase system, i.e. that both interphases can not collapse in such a way that a positive volume of fluid between the interphases it gets ejected in finite time. Lets assume that both surfaces collapse at time $T$ in a domain $D$ such that
$$\displaystyle{\lim_{t \to T-}} \,
[ f ( x_1 , x_2 , t ) - g ( x_1 , x_2 , t ) ] = 0\quad\text{for}\quad (x_1 - \widetilde{x}_1)^2 + (x_2 - \widetilde{x}_2)^2 < a^2
$$
where the constant $a>0$ and the point $\widetilde{x}=(\widetilde{x}_1,\widetilde{x}_2,\widetilde{x}_3)\in D$ are fixed.
Consider the domain $\Omega(t)$ denoted by
$$
\Omega(t) =\{(x_1,x_2,x_3): (x_1-\widetilde{x}_1)^2 + (x_2 - \widetilde{x}_2)^2 \leq ( R ( t ))^2 ,  \,
g ( x_1 , x_2, t) \leq x_3 \leq f ( x_1 , x_2 , t) \}
$$
with $0<R(t)< a$ and a section of its boundary given by
$$
{\cal{S}} ( t ) =\{(x_1,x_2,x_3): (x_1 - \widetilde{x}_1)^2 + (x_2-\widetilde{x}_2)^2 = ( R ( t ))^2 ,  \,
g ( x_1 , x_2, t) \leq x_3 \leq f (x_1,x_2,t)\}.
$$

Then by using the divergence free vector field $u$ it follows that
$$
\frac{d}{dt} \, {\rm Vol} \, \Omega ( t ) \, =  \,
\displaystyle{\int_{{\cal{S}}(t)}} \,
[ R^\prime ( t ) - u \cdot \nu ] \, d \; ( {\rm Area}) \,
$$
where $\nu$ is the unit normal to ${\cal{S}} ( t )$. If the integral $\int_0^T|u|_{L^{\infty}}dt$ is bounded then we can choose a time $t_0\in [0,T)$ and take
$$
R ( t ) = \frac{1}{2} \, a \,  - \,
\displaystyle{\int^{T}_{t}} \,
| u  |_{L^{\infty}} \, d \tau \  \ {\rm for} \  \ t_0 \leq t < T \,
$$
such that $0<R(t)<a$ for all $t\in [t_0,T)$. Consequently $\frac{d}{dt} \,
{\rm Vol} \, \Omega ( t ) \, \geq 0$ for all $t\in [t_0,T)$ which prevents a collapse (squirt singularity) forming
in between the moving surfaces (for more details see \cite{potato}). For a  general n-dimensional definition of a squirt singularity see \cite{squirts}. Therefore, the third part of the paper is devoted to show a bound of the velocity of the fluid in terms of $C^{1,\gamma}$ norms ($0<\gamma<1$) of the free boundary. The estimate is based on the property that in the principal value \eqref{sio} the mean of the kernels $K$ are zero on hemispheres. This extra cancellation was used by Bertozzi and Constantin \cite{B-C} for the vortex patch problem of the 2D Euler equation to prove no formation of singularities. For this system the convected vorticity takes constant values in disjoint domains and is related with the incompressible velocity by the Biot-Savart law (see \cite{M-B} for more details). Let us point out that the system we are dealing with is a more singular one. Also, we quote the work \cite{Verdera} of Mateu, Orobitg and  Verdera where this cancelation was used in quasiconformal mappings theory.

Finally we would like to emphasize how the character of the kernels becomes crucial here since for analogous active scalar models we can not obtain this result. For the 2D surface quasi-geostrophic  equation (SQG) \cite{CMT} the kernels are odd (Riesz transforms) and for the patch problem \cite{Y} the velocity is not in $L^\infty$, it is in $BMO$ (see \cite{St3} to get the definition and properties of the $BMO$ space). Furthermore in the case of regular initial data for SQG and the system \eqref{efm} the problem is also open \cite{CGO}.

The structure of the article is as follows. In section 2 we derived the equations of both interphases that are coupled together and in section 3 we show that the system is well-possed in Sobolev spaces. Finally, in section 4, we give a proof of boundedness of the velocity in terms of the smoothness of the interphases.

\section{The evolution equation of the moving surfaces}

The goal is to obtain the dynamics of a fluid that takes three different constant values of density $\rho^1,$ $\rho^2$ and $\rho^3$ as follows
\begin{equation}\label{frho}
 \rho(x_1,x_2,x_3,t)=\left\{\begin{array}{cl}
                    \rho^1\quad\mbox{in}&\{x_3>f(x,t)\},\\
                    \rho^2\quad\mbox{in}&\{f(x,t)>x_3>g(x,t)\},\\
                    \rho^3\quad\mbox{in}&\{g(x,t)>x_3\},

                 \end{array}\right.
\end{equation}
with $f(x,t)> g(x,t)$ for all $x=(x_1,x_2)\in \R^2$. Hence
$$
\grad \rho=(\rho^2-\rho^1)(\dxu f,\dxd f,-1)\delta(x_3-f(x,t))+(\rho^3-\rho^2)(\dxu g,\dxd g,-1)\delta(x_3-g(x,t)),
$$
where $\delta$ is the Dirac distribution defined by
\begin{equation*}
<h\delta,\eta>=\int_{\R^2} h(x)\eta(x,f(x,t))dx,
\end{equation*}
for $\eta(x_1,x_2,x_3)$ a test function.

For a divergence free velocity field, Darcy's law provides
\begin{equation}\label{sior}
u=(\dxu\Delta^{-1}\dxt\rho,\dxd\Delta^{-1}\dxt\rho,-\dxu\Delta^{-1}\dxu\rho-\dxd\Delta^{-1}\dxd\rho),
\end{equation}
and therefore
\begin{align}\label{pesao}
\begin{split}
u(x_1,x_2,x_3,t)=&-\D\frac{\rho^2-\rho^1}{4\pi}PV\int_{\R^2}\frac{(y_1,y_2,\grad f(x-y,t)\cdot y)}
{[|y|^2+(x_3-f(x-y,t))^2]^{3/2}}dy\\
&-\D\frac{\rho^2-\rho^1}{4\pi}PV\int_{\R^2}\frac{(y_1,y_2,\grad g(x-y,t)\cdot y)}
{[|y|^2+(x_3-g(x-y,t))^2]^{3/2}}dy,
\end{split}
\end{align}
for $f(x,t)\neq x_3\neq g(x,t)$. Taking $x_3\rightarrow f(x,t)$ yields
\begin{align*}
\begin{split}
u(x,f(x,t),t)=&-\D\frac{\rho^2-\rho^1}{4\pi}PV\int_{\R^2}\frac{(y_1,y_2,\grad f(x-y,t)\cdot y)}
{[|y|^2+(f(x,t)-f(x-y,t))^2]^{3/2}}dy\\
&-\D\frac{\rho^3-\rho^2}{4\pi}PV\int_{\R^2}\frac{(y_1,y_2,\grad g(x-y,t)\cdot y)}
{[|y|^2+(f(x,t)-g(x-y,t))^2]^{3/2}}dy,
\end{split}
\end{align*}
without considering tangential terms. In any case the evolution of the surfaces are given by the normal velocity, the tangential terms are related with the parameterization of the free boundary \cite{Hou}. Consequently
$$
(x,f(x,t))_t\cdot (-\dxu f,-\dxd f,1)=u(x,f(x,t),t)\cdot (-\dxu f,-\dxd f,1),
$$
and therefore
\begin{align}
\begin{split}\label{ce1}
f_t(x,t)=&\D\frac{\rho^2-\rho^1}{4\pi}PV\int_{\R^2}\frac{(\grad f(x,t)-\grad f(x-y,t))\cdot y)}
{[|y|^2+(f(x,t)-f(x-y,t))^2]^{3/2}}dy\\
&+\D\frac{\rho^3-\rho^2}{4\pi}PV\int_{\R^2}\frac{\grad f(x,t)\cdot y -\grad g(x-y,t)\cdot y}
{[|y|^2+(f(x,t)-g(x-y,t))^2]^{3/2}}dy.
\end{split}
\end{align}
In a similar way we obtain
\begin{align}
\begin{split}\label{ce2}
g_t(x,t)=&\D\frac{\rho^3-\rho^2}{4\pi}PV\int_{\R^2}\frac{(\grad g(x,t)-\grad g(x-y,t))\cdot y)}
{[|y|^2+(g(x,t)-g(x-y,t))^2]^{3/2}}dy\\
&+\D\frac{\rho^2-\rho^1}{4\pi}PV\int_{\R^2}\frac{\grad g(x,t)\cdot y -\grad f(x-y,t)\cdot y}
{[|y|^2+(g(x,t)-f(x-y,t))^2]^{3/2}}dy.
\end{split}
\end{align}
This coupled system describes the evolution of the moving boundaries $S_f (x,t ) $ and $S_g (x,t ) $.

\section{Well-posedness in the stable scenario $\rho^1 < \rho^2 < \rho^3$}

Let us define the function $d(f,g)(x,y,t)$ by the formula
$$
d(f,g)(x,y,t)=\frac{1}{[|y|^2+(f(x,t)-g(x-y,t))^2]^{1/2}}\qquad \forall\, x,y\in\R^2,
$$
which measure the distance between the contours $f$ and $g$. Therefore we consider $f(x,t)$ approaching the value $C_{\infty}$ as $|x|\rightarrow\infty$ to avoid that the surfaces collapse at infinity.

The section is devoted to prove the following theorem:
\begin{thm}\label{theorem}
Let $f_0(x)-C_\infty,g_0(x)\in H^{k}(\R^2)$ for $k\geq 4$, $d(f_0,g_0)\in L^\infty$ and $\rho_3>\rho_2>\rho_1$. Then there exists a time $T>0$ so
that there is a unique solution to \eqref{ce1} and \eqref{ce2} given by $(f(x,t),g(x,t))$ where $f(x,t)-C_\infty,g(x,t)\in C^{1}([0,T];H^k(\R^2))$, $f(x,0)=f_0(x)$ and $g(x,0)=g_0(x)$.
\end{thm}

Proof: We shall show that the following estimate holds
\begin{equation}\label{qtc}
\frac{d}{dt} E(t)\leq C(E(t)+1)^p
\end{equation}
for universal constants ($C$,$p$)  and the function $E(t)$ defined by
$$
E(t)=\|f-C_{\infty}\|_{H^k}(t)+\|g\|_{H^k}(t)+\|d(f,g)\|_{L^\infty}(t).
$$
Applying standard energy estimates argument permit us to conclude  local existence. This is based on introducing a regularized version of the system \eqref{ce1} and \eqref{ce2} that allows to take limits satisfying uniformly the a priori bound \eqref{qtc}.

 To simplify the exposition we shall consider $\rho^2-\rho^1=\rho^3-\rho^2=4\pi$ and $k=4$, being the rest of the cases analogous.  Some of the terms can be estimated exactly as in \cite{DY} and therefore we shall show below  how to deal with the different ones.

We have
\begin{align*}
\D\frac{1}{2}\dt\|f-C_\infty\|_{L^2}^2(t)&=\int_{\R^2}(f(x)-C_\infty) PV\int_{\R^2}\frac{(\grad f(x)-\grad f(x-y))\cdot
y}{[|y|^2+(f(x)-f(x-y))^2]^{3/2}}dydx\\
&\quad +\int_{\R^2}(f(x)-C_\infty)\grad f(x)\cdot PV\int_{\R^2}\frac{y}{[|y|^2+(f(x)-g(x-y))^2]^{3/2}}\\
&\quad -\int_{\R^2}(f(x)-C_\infty) PV\int_{\R^2}\frac{\grad g(x-y)\cdot y}{[|y|^2+(f(x)-g(x-y))^2]^{3/2}}\\
&=I_1+I_2+I_3.
\end{align*}
For $I_1$ we decompose further: $I_1=J_1+J_2+J_3$ where
\begin{align*}
J_1=\int_{|y|<1}\int_{\R^2}(f(x)-C_\infty)\frac{(\grad f(x)-\grad f(x-y))\cdot
y}{[|y|^2+(f(x)-f(x-y))^2]^{3/2}}dydx,
\end{align*}
\begin{align*}
J_2=\int_{|y|>1}(f(x)-C_\infty) PV\int_{|y|>1}\frac{\grad f(x)\cdot y}{[|y|^2+(f(x)-f(x-y))^2]^{3/2}}dydx
\end{align*}
and
\begin{align*}
J_3=\quad -\int_{y|>1}(f(x)-C_\infty) PV\int_{|y|>1}\frac{\grad f(x-y)\cdot y}{[|y|^2+(f(x)-f(x-y))^2]^{3/2}}dydx.
\end{align*}
Since $$\partial_{x_i}f(x)-\partial_{x_i}f(x-y)=\int_0^1\grad\partial_{x_i}f(x+(s-1)y)\cdot y\, ds,$$ we get
\begin{align*}
J_1&\leq  \int_0^1ds\int_{|y|<1}|y|^{-1}\int_{\R^2}\frac{|f(x)-C_{\infty}||\grad^2 f(x+(s-1)y)|}
{[1+((f(x)-f(x-y))^2|y|^{-2}]^{3/2}}dxdy\leq C \|f-C_\infty\|_{L^2}\|\grad f\|_{H^1}.
\end{align*}
In the term $J_2$ integration by parts yields
\begin{align*}
J_2&=\frac32\int_{|y|>1}\int_{\R^2}|f(x)-C_\infty|^2\frac{(f(x)-f(x-y))(\grad f(x)-\grad f(x-y))\cdot
y}{[|y|^2+((f(x)-f(x-y))^2]^{5/2}}dxdy\\
&\leq \frac32\int_{|y|>1}|y|^{-3}\int_{\R^2}|f(x)-C_\infty|^2\frac{|f(x)-f(x-y)||y|^{-1}(|\grad f(x)|+|\grad
f(x-y)|)}{[1+((f(x)-f(x-y))^2|y|^{-2}]^{5/2}}dxdy\\
&\leq C\|f-C_\infty\|^2_{L^2}\|\grad f\|^2_{L^\infty}\leq C\|f-C_\infty\|^2_{L^2}\|\grad f\|^2_{H^2}.
\end{align*}
In a similar way for $J_3$ we have
\begin{align*}
J_3&=\int_{|y|>1}\int_{\R^2}(f(x)-C_\infty)(f(x-y)-C_\infty)\frac{2|y|^2-(f(x)-f(x-y))^2}{[|y|^2+(f(x)-f(x-y))^2]^{5/2}}dxdy\\
&\quad+3\int_{|y|>1}\int_{\R^2}(f(x)-C_\infty)(f(x-y)-C_\infty)\frac{(f(x)-f(x-y))\grad f(x-y)\cdot y}{[|y|^2+(f(x)-f(x-y))^2]^{5/2}}dxdy\\
&\quad-\int_{|y|=1}\int_{\R^2}(f(x)-C_\infty)(f(x-y)-C_\infty)\frac{|y|^2}{[|y|^2+(f(x)-f(x-y))^2]^{3/2}}dxd\sigma(y)\\
&\leq C(\|\grad f\|_{L^{\infty}}+1)\|f-C_\infty\|^2_{L^2}.
\end{align*}
The term $I_2$ is estimated as follows
\begin{align*}
I_2&=\frac32\int_{\R^2}\int_{\R^2}(f(x)-C_\infty)^2\frac{(f(x)-g(x-y))(\grad f(x)-\grad g(x-y))\cdot y}{[|y|^2+(f(x)-g(x-y))^2]^{5/2}}dxdy\\
&=3\int_{|y|>1}\int_{\R^2} dxdy + 3\int_{|y|<1}\int_{\R^2} dxdy
\end{align*}
\begin{align*}
I_2&\leq \frac32\int_{|y|>1}|y|^{-3}\int_{\R^2}|f(x)-C_\infty|^2(|\grad f(x)|+|\grad g(x-y)|) dxdy\\
&\quad +\frac32\int_{|y|<1}\int_{\R^2}|f(x)-C_\infty|^2(|\grad f(x)|+|\grad g(x-y)|)(d(f,g)(x,y))^3 dxdy\\
&\leq C\|f-C_\infty\|^2_{L^2}(\|\grad f\|_{L^{\infty}}+\|\grad g\|_{L^{\infty}})(\|d(f,g)\|_{L^\infty}^3+1)
\end{align*}
For $I_3$ we write
\begin{align*}
I_3&=\int_{\R^2}\int_{\R^2}(f(x)-C_\infty)g(x-y)\frac{2|y|^2-(f(x)-g(x-y))^2}{[|y|^2+(f(x)-g(x-y))^2]^{5/2}}dxdy\\
&\quad +\int_{\R^2}\int_{\R^2}(f(x)-C_\infty)g(x-y)\frac{3(f(x)-g(x-y))(\grad f(x)-\grad g(x-y))\cdot y}{[|y|^2+(f(x)-g(x-y))^2]^{5/2}}dxdy
\end{align*}
and proceeding as before we find $$I_3\leq C\|f-C_\infty\|_{L^2}\|g\|_{L^2}(\|\grad f\|_{L^{\infty}}+\|\grad g\|_{L^{\infty}}+1)(\|d(f,g)\|_{L^\infty}^3+1).$$
Thus Sobolev inequalities yields
\begin{align}\label{ee1}
\D\dt\|f-C_\infty\|_{L^2}^2(t)&\leq C(E(t)+1)^{p+1}.
\end{align}
The analogous estimate for g is
\begin{align}\label{ee2}
\D\dt\|g\|_{L^2}^2(t)&\leq C(E(t)+1)^{p+1}.
\end{align}
To estimate the higher order derivative, we consider the quantity
\begin{align*}
\D\frac{1}{2}\dt\|\dxu^4f\|_{L^2}^2(t)&=I_4+I_5+I_6,
\end{align*}
where
\begin{align*}
I_4&=\int_{\R^2}\dxu^4f(x)\dxu^4\Big(PV\int_{\R^2}\frac{(\grad
f(x)-\grad f(x-y))\cdot y}{[|y|^2+(f(x)-f(x-y))^2]^{3/2}}dy\Big) dx,
\end{align*}
\begin{align*}
I_5&=\int_{\R^2}\dxu^4f(x)\dxu^4\Big(PV\int_{\R^2}\frac{\grad f(x)\cdot y}
{[|y|^2+(f(x)-g(x-y))^2]^{3/2}}dy\Big)dx
\end{align*}
and
\begin{align*}
I_6&=-\int_{\R^2}\dxu^4f(x)\dxu^4\Big(PV\int_{\R^2}\frac{\grad g(x-y)\cdot y}
{[|y|^2+(f(x,t)-g(x-y,t))^2]^{3/2}}dy \Big)dx.
\end{align*}
The estimates for $I_4$ are obtained  in \cite{DY}.
 
 We split $I_5$ and consider the most singular terms which are
  \begin{align*}
J_4&=\int_{\R^2}\dxu^4f(x)\grad \dxu^4f(x)\cdot PV\int_{\R^2}\frac{ y}{[|y|^2+(f(x)-g(x-y))^2]^{3/2}}dydx,
\end{align*}
\begin{align*}
J_5&=-3\int_{\R^2}\dxu^4f(x)\int_{\R^2}\frac{\grad f(x)\cdot y(f(x)-g(x-y))(\dxu^4f(x)-\dxu^4g(x-y))}
{[|y|^2+(f(x)-g(x-y))^2]^{5/2}}dydx,
\end{align*}
and
\begin{align*}
J_6&=C\int_{\R^2}\dxu^4f(x)\int_{\R^2}\frac{\grad f(x)\cdot y(f(x)-g(x-y))^4(\dxu f(x)-\dxu g(x-y))^4}
{[|y|^2+(f(x)-g(x-y))^2]^{11/2}}dydx.
\end{align*}
For integral $J_4$ we find that
\begin{align*}
J_4&=\frac32\int_{\R^2}\int_{\R^2}|\dxu^4 f(x)|^2\frac{(f(x)-g(x-y))(\grad f(x)-\grad g(x-y))\cdot y}{[|y|^2+(f(x)-g(x-y))^2]^{5/2}}dxdy\\
&\leq C\|\dxu^4 f\|_{L^2}^2(\|\grad f\|_{L^{\infty}}+\|\grad g\|_{L^{\infty}})(\|d(f,g)\|_{L^\infty}^3+1).
\end{align*}
For $J_5$ it follows
\begin{align*}
J_5&=\int_{|y|>1}|y|^{-3}\int_{\R^2}|\dxu^4f(x)||\grad f(x)|(|\dxu^4f(x)|+|\dxu^4g(x-y)|)dxdy\\
&\quad +\int_{|y|<1}\int_{\R^2}|\dxu^4f(x)||\grad f(x)|(|\dxu^4f(x)|+|\dxu^4g(x-y)|)(d(f,g)(x,y))^3dxdy\\
&\leq C\|\dxu^4 f\|_{L^2}\|\grad f\|_{L^{\infty}}(\|\dxu^4 f\|_{L^2}+\|\dxu^4 g\|_{L^{2}})(\|d(f,g)\|_{L^\infty}^3+1).
\end{align*}
Similarly we obtain
\begin{align*}
J_6&=C\int_{|y|>1}|y|^{-6}\int_{\R^2}|\dxu^4f(x)||\grad f(x)|(|\dxu f(x)|+|\dxu g(x-y)|)^4dxdy\\
&\quad +C\int_{|y|<1}\int_{\R^2}|\dxu^4f(x)||\grad f(x)|(|\dxu f(x)|+|\dxu g(x-y)|)^4(d(f,g)(x,y))^6dxdy\\
&\leq C\|\dxu^4 f\|_{L^2}\|\grad f\|_{L^{2}}(\|\dxu  f\|_{L^\infty}+\|\dxu g\|_{L^{\infty}})^4(\|d(f,g)\|_{L^\infty}^6+1)
\end{align*}
and finally
\begin{align*}
\D\dt\|\dxu^4 f\|_{L^2}^2(t)&\leq C(E(t)+1)^{p+1}.
\end{align*}
Analogously we get
\begin{align*}
\D\dt\|\dxd^4 f\|_{L^2}^2(t)&\leq C(E(t)+1)^{p+1},
\end{align*}
and therefore using \eqref{ee1} it follows
\begin{align}\label{esf}
\D\dt\|f-C_\infty\|_{H^4}^2(t)&\leq C(E(t)+1)^{p+1}.
\end{align}
In order to estimate the $H^4$ norm of $g$ we proceed as for the surface $f$ to obtain
\begin{align}\label{esg}
\D\dt\|g\|_{H^4}^2(t)&\leq C(E(t)+1)^{p+1}.
\end{align}
Next, we analyze the evolution of the distance between the surfaces. For the quantity $d(f,g)(x,y)$ we have
\begin{align}
\begin{split}\label{mct}
\D\dt d(f,g)(x,y,t)&=-\frac{(f(x)-g(x-y,t))(f_t(x)-g_t(x-y))}{[|y|^2+(f(x)-g(x-y))^2]^{3/2}}\\
&\leq (d(f,g)(x,y,t))^2(\|f_t\|_{L^\infty}(t)+\|g_t\|_{L^\infty}(t)).
\end{split}
\end{align}
We now estimate $\|f_t\|_{L^\infty}(t)$, being equivalent to control $\|g_t\|_{L^\infty}(t)$. We consider for $f_t(x,t)$ the following splitting:
$$
I_7=PV\int_{\R^2}\frac{(\grad f(x,t)-\grad f(x-y,t))\cdot y}
{[|y|^2+(f(x,t)-f(x-y,t))^2]^{3/2}}dy,
$$
$$
I_8=PV\int_{\R^2}\frac{\grad f(x,t)\cdot y}
{[|y|^2+(f(x,t)-g(x-y,t))^2]^{3/2}}dy,
$$
$$
I_9=-PV\int_{\R^2}\frac{\grad g(x-y,t)\cdot y}
{[|y|^2+(f(x,t)-g(x-y,t))^2]^{3/2}}dy.
$$
Let us decompose further: $I_7=J_7+J_8+J_9$ where
$$
J_7=\int_{|y|<1}\frac{(\grad f(x,t)-\grad f(x-y,t))\cdot y}
{[|y|^2+(f(x,t)-f(x-y,t))^2]^{3/2}}dy,
$$
$$
J_8=\grad f(x,t)\cdot PV\int_{|y|>1}\frac{ y}
{[|y|^2+(f(x,t)-f(x-y,t))^2]^{3/2}}dy
$$
and
$$
J_9=-PV\int_{|y|>1}\frac{\grad f(x-y,t)\cdot y}
{[|y|^2+(f(x,t)-f(x-y,t))^2]^{3/2}}dy.
$$
Hence
\begin{align*}
|J_7|&\leq C\int_{|y|<1}|y|^{-1}dy \|\grad^2  f\|_{L^\infty}\leq C \|\grad^2  f\|_{L^\infty}.
\end{align*}
We rewrite $J_8$ as
$$
J_8=\grad f(x,t)\cdot \int_{|y|>1}\frac{ y}{|y|^3}\frac{1-[1+(f(x,t)-f(x-y,t))^2|y|^{-2}]^{3/2}}
{[1+(f(x,t)-f(x-y,t))^2|y|^{-2}]^{3/2}}dy,
$$
and considering the function $Q(\al)=[1+\al^2]^{3/2}$, the mean value theorem gives
$$
J_8=\grad f(x,t)\cdot \int_{|y|>1}\frac{ y}{|y|^4}\frac{-3[1+\beta^2]^{1/2}\beta|f(x)-f(x-y)|}
{[1+(f(x,t)-f(x-y,t))^2|y|^{-2}]^{3/2}}dy,
$$
where $0\leq\beta\leq |f(x)-f(x-y)||y|^{-1}$. Therefore
$$
|J_8|\leq C\|\grad f\|_{L^\infty}\|f-C_\infty\|_{L^\infty}.
$$
In a similar manner as for $J_3$, integrations by parts in $J_9$ yields
$$
|J_9|\leq C(\|\grad f\|_{L^\infty}+1)\|f-C_\infty\|_{L^\infty}.
$$
For $I_8$ we split
\begin{align*}
I_8&=PV \int_{|y|>1} dy+\int_{|y|<1} dy =K_1+K_2.
\end{align*}
We deal with $K_1$ as we have done with $J_8$ to get $$|K_1|\leq \|\grad f\|_{L^\infty}(\|f-C_\infty\|_{L^\infty}+C_\infty+\|g\|_{L^\infty}).$$
For $K_2$ it is easy to obtain
$$|K_2|\leq \|\grad f\|_{L^\infty}\|d(f,g)\|^2_{L^\infty}.$$
We control $I_9$ as the term $I_3$ and therefore
$$
|I_9|\leq C\|g\|_{L^\infty}(\|\grad f\|_{L^{\infty}}+\|\grad g\|_{L^{\infty}}+1)(\|d(f,g)\|_{L^\infty}^3+1).
$$
The above estimates together with \eqref{mct} yields
\begin{align*}
\D\dt d(f,g)(x,y,t)&\leq d(f,g)(x,y,t)C(E(t)+1)^p,
\end{align*}
and integrating in time
\begin{align*}
d(f,g)(x,y,t+h)&\leq d(f,g)(x,y,t)\exp \Big(\int_t^{t+h}C(E(s)+1)^pds\Big),
\end{align*}
for $h>0$. Hence
\begin{align*}
\|d(f,g)\|_{L^\infty}(t+h)&\leq \|d(f,g)\|_{L^\infty}(t)\exp \Big(\int_t^{t+h}C(E(s)+1)^pds\Big).
\end{align*}
The above estimate applied to the following limit:
$$
\dt \|d(f,g)\|_{L^\infty}(t)=\lim_{h\rightarrow0^+}\frac{\|d(f,g)\|_{L^\infty}(t+h)-\|d(f,g)\|_{L^\infty}(t)}{h}
$$
allows us to get finally
\begin{align}
\begin{split}\label{mct2}
\D\dt \|d(f,g)\|_{L^\infty}(t)\leq C(E(t)+1)^p.
\end{split}
\end{align}
Combining estimates \eqref{esf}, \eqref{esg} and \eqref{mct2} we conclude that \eqref{qtc} holds for universal constants $C$ and $p$.

For uniqueness we proceed as in \cite{DY} which leads to the desired result.

\section{Bound for the fluid velocity: even kernels}

In this section we prove the following lemma:
\begin{lemma} Let $f,g$ be solutions of the contour system \eqref{ce1} and \eqref{ce2}. Then the velocity of the fluid satisfies the following bound
\begin{align}
\begin{split}\label{avel}
\|u\|_{L^\infty}\leq C\Big(1&+\frac{1}{\gamma}+\frac{1}{\gamma}\ln (1+\|\grad f\|_{C^\gamma}+\|\grad g\|_{C^\gamma})\\
&+\ln (1+\|f-C_{\infty}\|_{L^\infty}+C_\infty+\|\grad f\|_{L^2}
+\|g\|_{L^\infty}+\|\grad g\|_{L^2})\Big),
\end{split}
\end{align}
where $0<\gamma<1$ and the constant $C=C(\rho^1,\rho^2,\rho^3)$ depends on $\rho^1$, $\rho^2$ and $\rho^3$.
\end{lemma}

 Proof: We shall denote $x\in \R^2$ and $\overline{x}=(x,x_3)\in\R^3$. In order to acquire  the above inequality, we can split $\rho$ as follows:
$$
\rho=\rho^1\chi_{\Omega^+(t)}+\rho^1\chi_{\Omega^1(t)}+\rho^2\chi_{\Omega^2(t)}+(\rho^3-\rho^2)\chi_{\Omega^3(t)}+\rho^3\chi_{\Omega^-(t)}
$$
where $$\Omega^+(t)=\{x_3>M\},\qquad \Omega^-(t)=\{x_3<-M\},$$
$$\Omega^1(t)=\{M>x_3>f(x,t)\},\qquad\Omega^2(t)=\{f(x,t)>x_3>-M\}$$
and
$$\Omega^3(t)=\{g(x,t)>x_3>-M\}$$
for $M=\|f-C_{\infty}\|_{L^\infty}(t)+C_{\infty}+\|g\|_{L^\infty}(t)+1$.

From \eqref{sior} we check that
\begin{equation*}
(\dxt\Delta^{-1}\dxu,\dxt\Delta^{-1}\dxd,-(\dxu\Delta^{-1}\dxu+
\dxd\Delta^{-1}\dxd))(\chi_{\Omega^\pm(t)})=0.
\end{equation*}
Thus the Fourier transform yields
\begin{align}
\begin{split}\label{npares}
u&=T(\rho^1\chi_{\Omega^1(t)}+\rho^2\chi_{\Omega^2(t)}+(\rho^3-\rho^2)\chi_{\Omega^3(t)})\\
&\quad-\frac23(0,0,\rho^1\chi_{\Omega^1(t)}+\rho^2\chi_{\Omega^2(t)}+(\rho^3-\rho^2)\chi_{\Omega^3(t)}),
\end{split}
\end{align}
where
$$
T(h)(\overline{x})=\frac1{4\pi}PV\int_{\R^3}K(\overline{x}-\overline{y})h(\overline{y})d\overline{y},\qquad \overline{x}\in \R^3
$$
and
$$
K(\overline{x})=\left(3\frac{x_1x_3}{|\overline{x}|^5},3\frac{x_2x_3}{|\overline{x}|^5},\frac{2x_3^2-x_1^2-x_2^2}{|\overline{x}|^5}\right).
$$
Then we have to deal with
$$
\rho^1 T(\chi_{\Omega^1(t)}), \qquad \rho^2T(\chi_{\Omega^2(t)}) \quad\mbox{and} \qquad (\rho^3-\rho^2)T(\chi_{\Omega^3(t)}).
$$
We will consider $T(\chi_{\Omega^1(t)})$ since the other terms are analogous. We proceed as in \cite{B-C}. The key point is that the kernels in $T$ are even. Therefore in the principal value the mean of the kernels are zero on a hemisphere.
We consider the three different coordinates of $T(\chi_{\Omega^1(t)})=(T_1,T_2,T_3)$. Then for a fixed $\overline{x}\in\R^3$ we have
$$
T_1(\overline{x})=\frac3{4\pi}PV\int_{\Omega^1(t)}\frac{(x_1-y_1)(x_3-y_3)}{|\overline{x}-\overline{y}|^5}d\overline{y}.
$$
We take a distance $\delta$ given by $$\delta=\frac1{3(1+\|\grad f\|_{C^\gamma}+\|\grad g\|_{C^\gamma})^{1/\gamma}}.$$
Then if $d(\overline{x},\Omega^1(t))>\delta$ we consider $\Omega^1(t)=U^1(t)\cup U^2(t)$ for $$U^1(t)=\Omega^1\cap\{(x_1-y_1)^2+(x_2-y_2)^2\leq L^2\},$$
and $$U^2(t)=\Omega^1\cap\{(x_1-y_1)^2+(x_2-y_2)^2\geq L^2\},$$
where $L=2(1+\|f-C_\infty\|_{L^\infty}+C_\infty+\|\grad f\|_{L^2}+\|g\|_{L^\infty}+\|\grad g\|_{L^2})$.

The splitting $T_1=I_1+I_2$ for
$$I_1(\overline{x})=\frac3{4\pi}PV\int_{U^1(t)}\frac{(x_1-y_1)(x_3-y_3)}{|\overline{x}-\overline{y}|^5}d\overline{y},$$
$$I_2(\overline{x})=\frac3{4\pi}PV\int_{U^2(t)}\frac{(x_1-y_1)(x_3-y_3)}{|\overline{x}-\overline{y}|^5}d\overline{y}$$
gives
$$
|I_1(\overline{x})|\leq C\int_\delta^{\sqrt{2} L} r^{-1}d r\leq  C\ln (\sqrt{2} L/\delta).
$$
We write $I_2=\lim_{R\rightarrow +\infty}I_2^R$ such that
$$
I_2^R=\frac3{4\pi}\int_{U_R^2(t)} \partial_{y_1}(\frac{(x_3-y_3)}{|\overline{x}-\overline{y}|^3})d\overline{y},
$$
where $$U_R^2(t)=\Omega^1\cap\{L^2\leq (x_1-y_1)^2+(x_2-y_2)^2\leq R^2\}.$$ Then integration by parts gives
$$
I_2^R=\frac3{4\pi}\int_{\partial U_R^2(t)} \frac{(x_3-y_3)}{|\overline{x}-\overline{y}|^3} n_1 d\sigma(\overline{y}),
$$
and therefore
\begin{align*}
I_2^R &=\frac3{4\pi}\int_{L<|x-y|<R}\frac{\partial_{x_1}f(y)(x_3-f(y))dy}{(|x-y|^2+(x_3-f(y))^2)^{3/2}}
 -\frac3{4\pi}\int_0^{2\pi}\!\!\int^M_{f(x+Ly')}\frac{(x_3-y_3)L y'_1 d\theta dy_3}{(L^2+(x_3-y_3)^2)^{3/2}}\\
&\quad-\frac3{4\pi}\int_0^{2\pi}\!\!\int^M_{f(x+Ry')}\frac{(x_3-y_3)R y'_1 d\theta dy_3}{(R^2+(x_3-y_3)^2)^{3/2}}
\end{align*}
for $y'=(\cos\theta,\sin\theta)$. By the Cauchy-Schwarz inequality we find easily that
$$
|I_2^R|\leq \frac3{4\pi}\frac{\|\dxu f\|_{L^2}}{L}+\frac{3M}{L}+\frac{3M}{R}.
$$
We take now $R$ to infinity to get $|I_2|\leq C$ for $C$ a universal constant.

If $d(\overline{x},\Omega^1(t))<\delta$ we set $T_1=J_1+J_2$ where
$$J_1(\overline{x})=\frac3{4\pi}PV\int_{\Omega^1(t)\cap B_\delta(\overline{x})}\frac{(x_1-y_1)(x_3-y_3)}{|\overline{x}-\overline{y}|^5}d\overline{y}$$
and
$$J_2(\overline{x})=\frac3{4\pi}PV\int_{\Omega^1(t)\cap B^c_\delta(\overline{x})}\frac{(x_1-y_1)(x_3-y_3)}{|\overline{x}-\overline{y}|^5}d\overline{y}.$$
For $J_1$ one can proceed as in \cite{B-C}. It is clear that $J_1(\overline{x})= 0$ if $d(\overline{x},\partial\Omega^1(t))=d\geq\delta.$ Therefore we take $d<\delta$. We show the argument for the boundary of $\Omega^1(t)$ given by $x_3-f(x,t)=0$. The planar section trivializes. Then we define the following set of directions
$$
S_r(\overline{x})=\{\overline{z}\in \R^3\,:\,|\overline{z}|=1,\, \overline{x}+r\overline{z}\in \Omega^1(t),\, d\leq r<\delta\}
$$
Also, we pick the point $\overline{a}\in \R^3$ on the boundary, $a_3-f(a,t)=0$, such that $|\overline{a}-\overline{x}|=d$. We define the following hemisphere:
$$
S^+(\overline{x})=\{\overline{z}\in \R^3\,:\,|\overline{z}|=1,\, (-\grad f(a,t),1)\cdot \overline{z}\geq 0\}.
$$
For the directions $\overline{z}$ in the symmetric difference
$$D_r(\overline{x})=(S_r(\overline{x})\setminus S^+(\overline{x}))\cup (S^+(\overline{x})\setminus S_r(\overline{x})),$$
we define the angles $\varphi(\overline{z})$ by
$$
\sin\varphi(\overline{z})=\frac{(-\grad f(a,t),1)\cdot \overline{z} }{\sqrt{1+|\grad f(a,t)|^2}\,|\overline{z}|}.
$$
It is easy to check that $\sin\varphi(\overline{z})>0$ gives $x_3+rz_3-f(x+rz)<0$ and $\sin\varphi(\overline{z})<0$ gives $x_3+rz_3-f(x+rz)>0$.
Therefore
\begin{align*}
|\sin\varphi(\overline{z})|&=\Big|\frac{(-\grad f(a,t),1)\cdot (\overline{a}-\overline{x})}{\sqrt{1+|\grad f(a,t)|^2}\,r}+\frac{(-\grad f(a,t),1)\cdot (\overline{x}+r\overline{z}-\overline{a})}{\sqrt{1+|\grad f(a,t)|^2}\, r} \Big|\\
&\leq \frac{d}{r}+\Big|\frac{(-\grad f(a,t),1)\cdot (\overline{x}+r\overline{z}-\overline{a})}{\sqrt{1+|\grad f(a,t)|^2}\, r}+ \frac{(a_3-f(a,t)) -(x_3+rz_3-f(x+rz,t))}{\sqrt{1+|\grad f(a,t)|^2}\,r}\Big|\\
&\leq \frac{d}{r}+\frac{|f(x+rz,t)-f(a,t)-\grad f(a)\cdot(x+rz-a)|}{r}
\end{align*}
and finally
\begin{align*}
|\sin\varphi(\overline{z})|&\leq \frac{d}{r}+\frac{\|\grad f\|_{C^\gamma}|x+rz-a|^{1+\gamma}}{r}\leq \frac{d}{r}+\frac{\|\grad f\|_{C^\gamma}(d+r)^{1+\gamma}}{r}\\
&\leq \frac{d}{r}+2^\gamma\frac{\|\grad f\|_{C^\gamma}(d^{1+\gamma}+r^{1+\gamma})}{r}\leq \frac{d}{r}+2^\gamma\frac{d^{1+\gamma}+r^{1+\gamma}}{\delta ^\gamma r}\leq (1+2^\gamma)\frac{d}{r}+\Big(\frac{2r}{\delta}\Big)^\gamma.
\end{align*}
The above estimate and the fact that
$$
\int_{S^+(\overline{x})}\frac{(x_1-y_1)(x_3-y_3)}{|\overline{x}-\overline{y}|^5}d\sigma(\overline{y})=0
$$
yield for $J_1$ the following bound
$$
|J_1(\overline{x})|\leq \frac{3}{4\pi}\int_d^\delta\int_0^{2\pi}\int_{-\pi/2}^{\pi/2} \chi_{\{\varphi(\overline{z}):\,\overline{z}\in D_r(\overline{x})\}}|\cos \varphi|d\varphi d\theta \frac{dr}{r}\leq \frac{3\pi}{2}(1+2^\gamma(1+1/\gamma)).
$$

We write $J_2=K_1+K_2$ where
$$K_1(\overline{x})=\frac3{4\pi}PV\int_{U^1(t)\cap B^c_\delta(\overline{x})}\frac{(x_1-y_1)(x_3-y_3)}{|\overline{x}-\overline{y}|^5}d\overline{y},$$
and
$$K_2(\overline{x})=\frac3{4\pi}PV\int_{U^2(t)\cap B^c_\delta(\overline{x})}\frac{(x_1-y_1)(x_3-y_3)}{|\overline{x}-\overline{y}|^5}d\overline{y}.$$
Then one can deal as for $I_1$ and $I_2$ respectively to obtain the same estimates.

We write the second coordinate as follows
$$
T_2(\overline{x})=\frac3{4\pi}PV\int_{\Omega^1(t)}\partial_{y_2}(\frac{(x_3-y_3)}{|\overline{x}-\overline{y}|^3})d\overline{y},
$$
to get an analogous bounds.

For the third coordinate we have
$$
T_3(\overline{x})=\frac1{4\pi}PV\int_{\Omega^1(t)}\frac{2(x_3-y_3)^2-(x_1-y_1)^2-(x_2-y_2)^2}{|\overline{x}-\overline{y}|^5}d\overline{y},
$$
and we can proceed as before. But for the term $I_2^R$  we have

\begin{align*}
I_2^R&=\frac1{4\pi}\int_{U_R^2(t)} \partial_{y_3}(\frac{(x_3-y_3)}{|\overline{x}-\overline{y}|^3})d\overline{y}=\frac1{4\pi}\int_{L<|x-y|<R}\int_{f(y)}^M \partial_{y_3}(\frac{(x_3-y_3)}{|\overline{x}-\overline{y}|^3}) dy_3dy\\
&=\frac1{4\pi}\int_{L<|x-y|<R}(P(\frac{x_3-M}{|x-y|})-P(\frac{x_3-f(y)}{|x-y|})) \frac{dy}{|x-y|^2},
\end{align*}
where $P(\al)=\al/(1+\al^2)^{3/2}$. The mean value theorem yields
$$
I_2^R=\frac1{4\pi}\int_{L<|x-y|<R}P'(\widetilde{\al}) \frac{(f(y)-M)dy}{|x-y|^3}
$$
and using that $|P'(\al)|\leq 2$ we find easily $|I_2^R|\leq 2M/L\leq 1$.

\subsection*{{\bf Acknowledgments}}

\smallskip

The authors were partially supported by the grant {\sc MTM2008-03754} of the MCINN (Spain) and
the grant StG-203138CDSIF  of the ERC. The second author was partially supported by NSF-DMS grant
0901810.

\begin{quote}
\begin{tabular}{ll}
Diego C\'ordoba &  Francisco Gancedo\\
{\small Instituto de Ciencias Matem\'aticas} & {\small Department of Mathematics}\\
{\small Consejo Superior de Investigaciones Cient\'ificas} & {\small University of Chicago}\\
{\small Serrano 123, 28006 Madrid, Spain} & {\small 5734 University Avenue, Chicago, IL 60637}\\
{\small Email: dcg@icmat.es} & {\small Email: fgancedo@math.uchicago.edu}
\end{tabular}
\end{quote}



\begin{thebibliography}{99}

\bibitem{bear} J. Bear, Dynamics of Fluids in Porous Media, \emph{American
Elsevier}, New York, 1972.

\bibitem{B-C} A.~L. Bertozzi and P. Constantin. Global regularity for vortex patches.
\emph{Comm. Math. Phys.} 152 (1): 19--28, 1993.

\bibitem{CMT} P.~Constantin, A.~J. Majda, and E.~Tabak. \newblock Formation
of strong fronts in the 2-{D} quasigeostrophic thermal active scalar.
\newblock \emph{Nonlinearity}, 7:1495--1533, 1994.

\bibitem{Peter} P. Constantin and M. Pugh. Global solutions for small data to the
Hele-Shaw problem. \emph{Nonlinearity}, 6 (1993), 393 - 415.

\bibitem{potato}D. C\'{o}rdoba and C. Fefferman. Potato chip
singularities of 3D flows. SIAM J. Math. Anal. 33 (2001)
no.4, 786-789

\bibitem{squirts}  D.~C\'ordoba, C.~Fefferman and R. de la Llave.
\newblock On squirt singularities in hydrodynamics.
 SIAM J. Math. Anal.  36  (2004),  no. 1, 204--213 .

\bibitem{DY} D. C\'ordoba and F. Gancedo. Contour dynamics of incompressible 3-D fluids in a porous medium with
different densities. \emph{Comm. Math. Phys},.  273  (2007),  no. 2, 445--471.

\bibitem{DY2} D. C\'ordoba and F. Gancedo. A maximum principle for the Muskat problem for fluids with different densities. \emph{Comm. Math. Phys.}, 286 (2009), no. 2, 681-696

\bibitem{CGO} D. C\'ordoba, F. Gancedo and R. Orive. Analytical behaviour of 2D incompressible flow in porous media. \newblock\emph{J. Math. Phys.}, 48, 6, 2007.

\bibitem{DPS} T. Dombre, A. Pumir and E. Siggia. On the interface dynamics for convection in porous media. \newblock\emph{Physica D}, 57 (1992) 311-329.

\bibitem{ES} J. Escher and G. Simonett. Classical solutions for Hele-Shaw models
with surface tension. \emph{Adv. Differential Equations}, (1997) 2:619-642.

\bibitem{Y} F. Gancedo. Existence for the $\alpha$-patch model and the QG sharp front in Sobolev spaces. \emph{Adv. Math.}, Vol 217/6: 2569-2598, 2008.

\bibitem{Hou} T.Y. Hou, J.S. Lowengrub and M.J. Shelley. Removing the Stiffness
from Interfacial Flows with Surface Tension. \emph{J. Comput. Phys.}, 114: 312-338, 1994.

\bibitem{M-B} A.~J. Majda and A.~L. Bertozzi. \newblock Vorticity
and the Mathematical Theory of Incompressible Fluid Flow. \newblock \emph{Cambridge Press}, 2002.

\bibitem{Verdera}  J. Mateu, J. Orobitg and J. Verdera \newblock Extra cancellation of even Calderon-Zygmund operators and quasiconformal mappings. \newblock \emph{arXiv:0802.1185v2 .}

\bibitem{Muskat} M. Muskat. \newblock The flow of homogeneous fluids through porous media. \newblock \emph{New York}, 1937.

\bibitem{S-T} P.G. Saffman and Taylor. The penetration of a fluid into a porous medium or Hele-Shaw cell containing a more viscous liquid. \newblock \emph{Proc. R. Soc. London, Ser. A} 245, 312-329, 1958.

\bibitem{SCH} M. Siegel, R. Caflisch and S. Howison. Global
Existence, Singular Solutions, and Ill-Posedness for the Muskat Problem. \emph{Comm. Pure and Appl.
Math.}, 57: 1374-1411, 2004.

\bibitem{St3} E.~Stein. \newblock Harmonic Analysis. \newblock \emph{%
Princeton University Press.} Princeton, NJ, 1993.


\end{thebibliography}
\end{document}